\documentclass[twoside,10pt]{article}

\usepackage{psfrag}
\usepackage{graphicx}
\usepackage{amssymb}
\usepackage{latexsym}

\newtheorem{theorem}{Theorem}
\newtheorem{proposition}{Proposition}
\newtheorem{definition}{Definition}
\newtheorem{lemma}{Lemma} 
\newtheorem{corollary}{Corollary}
\newtheorem{sublemma}{Sublemma}

\begin{document}

\title{Minkowski- Versus Euclidean Rank for Products of Metric Spaces }
\maketitle

\vspace{0.2cm}

\begin{center}
{\large Thomas Foertsch* $\hspace{1cm}$ and $\hspace{1cm}$ Viktor Schroeder}
\footnote{* supported by SNF Grant 21 - 589 38.99}
\end{center}

\vspace{0.5cm}


\begin{abstract}
We introduce a notion of the Euclidean- and the Minkowski rank for arbitrary metric spaces and we study
their behaviour with respect to products.
We show that the Minkowski rank is additive with respect to metric products, while additivity of the 
Euclidean rank only holds under additional assumptions, e.g. for Riemannian manifolds.
We also study products with nonstandard product metrics.
\end{abstract}

\vspace{0.2cm}


\section{Introduction}

For Riemannian manifolds there are various definitions of a rank in the literature
(compare e.g. \cite{bbe}, \cite{ss}, \cite{g}). A notion which can easily be generalized to arbitrary 
metric spaces is the rank as the maximal dimension of an Euclidean subspace isometrically embedded 
into the manifold.\\
It is known (an will also be proved below) that for Riemannian manifolds this {\it Euclidean rank}
is additive with respect to products. This is not the case for more general metric spaces, even 
for Finsler manifolds (see Theorem \ref{theoeukl} below). \\
In contrary it turns out that the {\it Minkowski rank} defined as the maximal dimension of an 
isometrically embedded normed vector space has a better functional behaviour with respect
to metric products.
\begin{definition} Minkowski- and Euclidean rank for metric spaces
\begin{description}
\item[a)] For an arbitrary metric space $(X,d)$ the  {\bf Minkowski rank} is
\begin{displaymath}
rank_{M}(X,d) \; := \; 
\sup_{(V,||\cdot ||)} \Big\{ dim V \; \Big| \; \exists \; \mbox{isometric map} \;\; 
i_V:(V,||\cdot ||) \longrightarrow (X,d) \Big\} .
\end{displaymath}
\item[b)] The {\bf Euclidean rank} is defined as
\begin{displaymath}
rank_{E}(X,d) \; := \; 
\sup \Big\{ n \in \mathbb{N} \; \Big| \; \exists \; \mbox{isometric map} \;\; 
i_{\mathbb{E}^n}:\mathbb{E}^n\longrightarrow (X,d) \Big\} .
\end{displaymath}
\end{description}
\end{definition} 

In the special case of Riemannian manifolds these rank definitions coincide.
\begin{theorem}  \label{theoriem}
Let ${\cal M}$ be a Riemannian manifold, then 
\begin{displaymath} 
rank_{M}({\cal M}) \; = \; rank_{E}({\cal M}).
\end{displaymath}
\end{theorem}
For more general metric spaces, the ranks may be different and they even have different functional
behaviour with respect to metric products. \\
The Minkowski rank is additive, i.e., we have
\begin{theorem} \label{theomink}
Let $(X_i,d_i)$, $i=1,2$, be metric spaces and denote their metric product by $(X_1\times X_2,d)$.
Then 
\begin{displaymath}
rank_M(X_1,d_1) \; + \; rank_M(X_2,d_2) \; = \; rank_M(X_1\times X_2,d).
\end{displaymath}
\end{theorem}
As an immediate consequence of Theorem \ref{theoriem} and Theorem \ref{theomink} we obtain the 
additivity of the Euclidean rank for Riemannian manifolds. 
\begin{corollary} \label{corriem}
Let $({\cal M}_i,g_i)$, $i=1,2$, be Riemannian manifolds and denote their Riemannian product 
by $({\cal M}_1\times {\cal M}_2,g)$. Then it holds
\begin{displaymath}
rank_E({\cal M}_1,g_1) \; + \; rank_E({\cal M}_2,g_2) \; = \; rank_E({\cal M}_1\times {\cal M}_2,g).
\end{displaymath}
\end{corollary}
In the general case the additivity of the Euclidean rank does not hold. In section
\ref{chapter-euclid} we give an 
example of two normed vector spaces $(V_i,||\cdot ||_i)$, $i=1,2$,
that do not admit an isometric embedding of $\mathbb{E}^2$, although $\mathbb{E}^3$ may be embedded in
their product. Thus $rank_E(V_i)=1$ for $i=1,2$ but $rank_E (V_1\times V_2)\ge 3$ and we 
obtain:
\begin{theorem} \label{theoeukl}
Let $(X_i,d_i)$, $i=1,2$, be metric spaces and denote their metric product by $(X_1\times X_2,d)$.
Then it holds
\begin{displaymath}
rank_E(X_1,d_1) \; + \; rank_E(X_2,d_2) \; \le \; rank_E(X_1\times X_2,d),
\end{displaymath}
but there are examples such that the inequality is strict.
\end{theorem}
In the results above the metric $d$ on the product $X_1\times X_2$ is the standard one.
More generally, let $(X_i,d_i)$, $i=1,...,n$, be metric spaces then there are different 
possibilities to define a metric $d$ on the product $X=\Pi_{i=1}^nX_i$. It is natural to
require that the metric on $X$ is of the form $d=d_{\Phi}$,
\begin{displaymath}
d_{\Phi} \Big( (x_1,...,x_n),(y_1,...,y_n) \Big) \; = \; 
\Phi \Big( d_1(x_1,y_1),...,d_n(x_n,y_n) \Big) ,
\end{displaymath}
where $\Phi :Q^n\longrightarrow [0,\infty )$ is a function defined on the quadrant 
$Q^n=[0,\infty )^n$. \\
If we require in addition that $(X,d_\Phi )$ is an inner metric space as long as all 
factors $(X_i,d_i)$ are inner metric spaces, then $\Phi$ has to be of a very special
type. We discuss this in detailed form in section \ref{product-metrics} and obtain in
particular
\begin{theorem} \label{innermetric}
Let $\Phi : Q^n\longrightarrow [0,\infty )$ be a function. Then $(X,d_{\Phi})$ is
an inner metric space for all possible choices of inner metric spaces $(X_i,d_i)$
if and only if 
\begin{displaymath}
\Psi : \mathbb{R}^n \longrightarrow [o,\infty ), \hspace{1,5cm}
\Psi \Big( \sum\limits_{i=1}^n \; x_i \, e_i\Big) \; := \; 
\Phi \Big( \sum\limits_{i=1}^n \; |x_i| \, e_i \Big)
\end{displaymath}
is a norm.
\end{theorem}
We prove the additivity of the Minkowski rank with respect to these generalized products:
\begin{theorem} \label{generalized-minkadd}
Let $\Phi : Q^n\longrightarrow [0,\infty )$ be a function such that
$\Psi : \mathbb{R}^n \longrightarrow [0,\infty )$ defined as in Theorem \ref{innermetric}
is a norm with a strictly convex norm ball. Let $(X_i,d_i)$, $i=1,...,n$, be metric spaces
and $X={\Pi}_{i=1}^nX_i$. Then
\begin{displaymath}
rank_{M}\Big( X,d_{\Phi}\Big) \; = \; \sum\limits_{i=1}^n \; rank_{M} \Big( X_i,d_i\Big) .
\end{displaymath} 
\end{theorem}

Finally we want to thank Andreas Bernig for useful discussions.


\section{Minkowski Rank for Riemannian Manifolds}

In this section we give a \\
{\bf Proof of Theorem \ref{theoriem}:}\\
Let $({\cal M},g)$ be a connected Riemannian manifold with induced distance function
$d:{\cal M}\times {\cal M} \longrightarrow [0,\infty )$. 
Obviously $rank_E({\cal M},d)\le rank_M({\cal M},d)$ and it remains to prove the opposite
inequality. \\
Let $k=rank_M({\cal M},d)$ and thus there exists an isometric map $\varphi : V\longrightarrow {\cal M}$,
where $(V,||\cdot ||)$ is a $k$-dimensional normed vector space. Clearly 
$W:=\varphi (V) \subset {\cal M}$ with the induced topology is homeomorphic to $\mathbb{R}^k$. \\
We first show that $W$ is a convex subset of ${\cal  M}$. Let therefore 
$\varphi (v_1), \varphi (v_2) \in W$ and consider the curve $c:\mathbb{R}\longrightarrow {\cal M}$,
$c(t)=\varphi(tv_2+(1-t)v_1)$. Since $t\longrightarrow tv_2+(1-t)v_2$ is a minimal geodesic in the space
$(V,||\cdot ||)$ and $\varphi$ is isometric, $c$ is a minimal geodesic in ${\cal M}$ and in particular the restriction
$c|_{[0,1]}$ is the (up to parametrization) unique minimal geodesic from $\varphi (v_1)$ to 
$\varphi (v_2)$ and contained in $W$.\\
By Theorem 1.6. in \cite{cg} $W$ is a totally geodesic submanifold which is in 
addition homeomorphic to $\mathbb{R}^k$. In particular $W$ is itself a Riemannian manifold with 
the induced metric, and $\varphi : V\longrightarrow W$ is an isometry. \\
Note that the abelian group $V$ acts on $W$ transitively by isometries via the action 
$\psi : V\times W \longrightarrow W$
\begin{displaymath}
\psi \Big( v_1,\varphi (v_2) \Big) \; = \; \varphi \Big( v_1+v_2 \Big) .
\end{displaymath}
Thus $W$ is a homogeneous Riemannian manifold, homeomorphic to $\mathbb{R}^k$ with a transitively
acting abelian group of isometries. Thus $W$ is isometric to the Euclidean space $\mathbb{E}^k$
and hence the image of an isometric map $\chi : \mathbb{E}^k \longrightarrow {\cal M}$.
Thus $rank_E({\cal M})\ge k$.
\begin{flushright}
{\bf q.e.d.}
\end{flushright}


\section{Minkowski Rank of Products I}

\label{minkadd}

In this section we prove that the Minkowski rank is additive for metric products.
Let therefore $(X_i,d_i)$, $i=1,2$, be metric spaces and consider the product $X=X_1\times X_2$
with the standard product metric 
\begin{displaymath}
d\Big( (x_1,x_2),(x_1',x_2') \Big) \; = \; 
\Big( d_1^2(x_1,x_1') \; + \; d_2^2(x_2,x_2'){\Big)}^{\frac{1}{2}}.
\end{displaymath}
We need an auxiliary result: Let $V$ be a real vector space and denote by $A$ the affine space
on which $V$ acts simply transitively. Thus for $a\in A$ and $v\in V$ the point $a+v\in A$ and 
for $a,b\in A$ the vector $b-a\in V$ are defined. As usual a pseudonorm on $V$ is a function 
$||\cdot ||$ which satisfies the properties of a norm with the possible exception that 
$||v||=0$ does not necessarily imply $v=0$. A pseudonorm $||\cdot ||$ on $V$ induces a pseudometric 
$d$ on $A$ via
\begin{displaymath}
d(a,b) \; = \; ||b-a|| \hspace{2cm} \forall \; a,b \in A.
\end{displaymath}

We denote the resulting pseudometric space by $(A,||\cdot ||)$. With this notation
we have:
\begin{proposition} \label{prop}
Let $(X_i,d_i)$, $i=1,2$, be metric spaces and
$\varphi : A \longrightarrow X_1 \times X_2$, $\varphi =({\varphi}_1,{\varphi}_2)$ be an 
isometric map. Then there exist pseudonorms $||\cdot ||_i$, $i=1,2$ on $V$, such that 
\begin{description}
\item[i)] $||v||^2 \; = \; ||v||_1^2 \; + \; ||v||_2^2 \hspace{1.5cm}$ and 
\item[ii)] ${\varphi}_i : (A,||\cdot ||_i) \longrightarrow (X_i,d_i)$, $i=1,2$ are isometric.
\end{description}
\end{proposition}
For the proof of Proposition \ref{prop} we define 
${\alpha}_i:A\times V \longrightarrow [0,\infty )$, $i=1,2$, via
\begin{displaymath}
{\alpha}_i(a,v) \; := \; d_i\Big( {\varphi}_i(a), {\varphi}_i(a+v) \Big) .
\end{displaymath}
Since $\varphi$ is isometric, we have
\begin{equation} \label{quadratgleichung}
{\alpha}_1^2(a,v) \; + \; {\alpha}_2^2(a,v) \; = \; d^2 \Big( \varphi (a), \varphi (a,v) \Big) \; = \; ||v||^2.
\end{equation}

We will prove the following Lemmata:

\begin{lemma} \label{lemma1}
\begin{displaymath}
{\alpha}_i(a,v) \; = \; {\alpha}_i(a+v,v), \;\; i=1,2, \hspace{2cm} \forall a\in A, v\in V,
\end{displaymath}
\end{lemma}

\begin{lemma} \label{lemma2}
\begin{displaymath}
{\alpha}_i(a,tv) \; = \; |t| {\alpha}_i(a,v), \;\; i=1,2, \hspace{2cm} \forall a\in A, v\in V, t\in \mathbb{R},
\end{displaymath} 
\end{lemma}

\begin{lemma} \label{lemma3}
\begin{displaymath}
{\alpha}_i(a,v) \; = \; {\alpha}_i(b,v), \;\; i=1,2, \hspace{2cm} \forall a,b\in A, v\in V \;\; \mbox{and}
\end{displaymath}
\end{lemma}

\begin{lemma} \label{lemma4}
\begin{displaymath}
{\alpha}_i(v+w) \; \le \; {\alpha}_i(v) \; + {\alpha}_i(w), \;\; i=1,2, \hspace{2cm} \forall v,w\in V,
\end{displaymath}
\end{lemma}
where ${\alpha}_i(v):={\alpha}_i(a,v)$ with $a\in A$ arbitrary (compare with Lemma \ref{lemma3}). \\

From Lemmata \ref{lemma1} - \ref{lemma4} it follows immediately, that $||\cdot||_i$ defined via
$||v ||_i:={\alpha}_i(v)$ $\forall v\in V$, $i=1,2$, is a pseudonorm on V. Furthermore from
\begin{eqnarray*}
d_i \Big( {\varphi}_i(a),{\varphi}_i(b) \Big) & = & 
d_i \Big( {\varphi}_i(a),{\varphi}_i(a+(b-a)) \Big) \\
& = & {\alpha}_i(b-a) \\
& = & ||b-a||_i \hspace{1cm} \forall a,b\in A
\end{eqnarray*}
it follows that 
\begin{displaymath}
{\varphi}_i : \; \Big( A,||\cdot ||_i\Big) \; \longrightarrow \; \Big( X_i,||\cdot ||_i\Big),
\hspace{1cm} i=1,2,
\end{displaymath}
are isometric mappings.\\

{\bf Proof of Lemma \ref{lemma1}:} \\
The $d_i$'s triangle inequality yields
\begin{equation}
{\alpha}_i(a,v) \; + \; {\alpha}_i(a+v,v) \; \ge \; {\alpha}_i(a,2v)
\end{equation}
and thus
\begin{equation} \label{quadr}
{\alpha}_i^2(a,v) \; + \; 2 {\alpha}_i(a,v) {\alpha}_i(a+v,v) \; + \; {\alpha}_i^2(a+v,v) \; \ge \;
{\alpha}_i^2(a,2v). 
\end{equation}
Using equation (\ref{quadratgleichung}) the sum of the equations (\ref{quadr}) for $i=1$ and $i=2$ 
becomes
\begin{displaymath}
||v||^2 \; + \; 2 \; \big<
\left(
\begin{array}{c}
{\alpha}_1(a,v) \\ {\alpha}_2(a,v)
\end{array}
\right),
\left(
\begin{array}{c}
{\alpha}_1(a+v,v) \\ {\alpha}_2(a+v,v)
\end{array}
\right)
\big> \; + \; ||v||^2 \; \ge \; 4 ||v||^2,
\end{displaymath}  
where $<\cdot ,\cdot >$ denotes the standard scalar product on ${\mathbb{R}}^2$.
Thus we have
\begin{displaymath}
\big<
\left(
\begin{array}{c}
{\alpha}_1(a,v) \\ {\alpha}_2(a,v)
\end{array}
\right),
\left(
\begin{array}{c}
{\alpha}_1(a+v,v) \\ {\alpha}_2(a+v,v)
\end{array}
\right)
\big> 
\; \ge \; ||v||^2.
\end{displaymath}

The Euclidean norm of the vectors $({\alpha}_1(a,v),{\alpha}_2(a,v))$ and
$({\alpha}_1(a+v,v),{\alpha}_2(a+v,v))$ equals $||v||$, due to equation \ref{quadratgleichung}.
Therefore the Cauchy Schwarz inequality yields 
\begin{displaymath}
\left(
\begin{array}{c}
{\alpha}_1(a,v) \\ {\alpha}_2(a,v)
\end{array}
\right)
\; = \;
\left(
\begin{array}{c}
{\alpha}_1(a+v,v) \\ {\alpha}_2(a+v,v)
\end{array}
\right).
\end{displaymath}
\begin{flushright}
{\bf q.e.d.}
\end{flushright}

{\bf Proof of Lemma \ref{lemma2}:}\\
The $d_i$'s triangle inequality yields for all $n\in \mathbb{N}$
\begin{displaymath}
{\alpha}_i(a,nv) \; \le \; \sum\limits_{k=0}^{n-1}{\alpha}_i(a+kv,v) \; = \; n{\alpha}_i (a,v),
\end{displaymath}
where the last equation follows from Lemma \ref{lemma1} by induction. Thus we find
\begin{eqnarray*}
n^2||v||^2 \; = \; ||nv||^2 & = &
{\alpha}_1^2(a,nv) \; + \; {\alpha}_2^2(a,nv) \\
& \le & n^2\Big( {\alpha}_1^2(a,v) \; + \; {\alpha}_2^2(a,v) \Big) \\
& = & n^2||v||^2 \hspace{1cm} \forall n\in \mathbb{N}, v\in V, a\in A
\end{eqnarray*}
and therefore
\begin{displaymath}
{\alpha}_i(a,nv) \; = \; n{\alpha}_i (a,v), \;\; i=1,2, \hspace{2cm} \forall n\in \mathbb{N}, v\in V,
a\in A. 
\end{displaymath}
Thus for $p,q \in \mathbb{N}$, it is 
\begin{displaymath}
q{\alpha}_i(a,\frac{p}{q}v) \; = \; {\alpha}_i(a,pv) \; = \; p{\alpha}_i(a,v),
\end{displaymath}
i.e.
\begin{displaymath}
{\alpha}_i(a,tv) \; = \; t{\alpha}_i(a,v) \hspace{2cm} \forall t\in \mathbb{Q}_+
\end{displaymath}
and by continuity even $\forall t \in \mathbb{R}_+$. \\
Finally note that for all $t\in \mathbb{R}_+$
\begin{displaymath}
{\alpha}_i(a,-tv) \; = \; {\alpha}_i(a-tv,tv) \; = \; {\alpha}_i(a,tv) \; = \; t{\alpha}_i(a,v), \;\; 
i=1,2,
\end{displaymath}
where the first equality is just the symmetry of the metric $d_i$ and the second equality follows 
from Lemma \ref{lemma1}.
\begin{flushright}
{\bf q.e.d.}
\end{flushright}

{\bf Proof of Lemma \ref{lemma3}:} \\
For $n\in \mathbb{N}$ we have
\begin{eqnarray*}
\Big|{\alpha}_i(a,nv) \; - \; {\alpha}_i(b,nv)\Big| & = & 
\Big| d_i\Big( {\varphi}_i(a),{\varphi}(a+nv)\Big) - d_i\Big( {\varphi}_i(b),{\varphi}(b+nv)\Big) \Big| \\
& \le & d_i \Big( {\varphi}_i(a),{\varphi}_i(b)\Big) \; + \;  
d_i \Big( {\varphi}_i(a+nv),{\varphi}_i(b+nv)\Big) \\
& \le &  d \Big( {\varphi}(a),{\varphi}(b)\Big) \; + \; d \Big( {\varphi}(a+nv),{\varphi}(b+nv)\Big) \\
& = & 2||b-a||, \;\;\;\;\;\;\; i=1,2,
\end{eqnarray*}
and therefore
\begin{displaymath}
{\alpha}_i(a,v) \; = \; \lim\limits_{n\longrightarrow \infty} \frac{1}{n}{\alpha}_i(a,nv)
\; = \;  \lim\limits_{n\longrightarrow \infty} \frac{1}{n}{\alpha}_i(b,nv) \; = \; 
{\alpha}_i(b,v), \;\; i=1,2.
\end{displaymath}
\begin{flushright}
{\bf q.e.d.}
\end{flushright}

{\bf Proof of Lemma \ref{lemma4}:} \\
The claim simply follows by
\begin{displaymath}
{\alpha}_i(v+w) \; = \; {\alpha}_i(a,v+w) \; \le \; {\alpha}_i(a,v) \; + \; {\alpha}_i(a+v,w) \; = \; 
{\alpha}_i(v) \; + \; {\alpha}_i(w),
\end{displaymath}
where the inequality follows by the $d_i$'s triangle inequality and the last equation is due to
Lemma \ref{lemma3}.
\begin{flushright}
{\bf q.e.d.}
\end{flushright}

With that we are now ready for the \\

{\bf Proof of Theorem \ref{theomink}:} 
\begin{description}
\item[i)] Superadditivity follows as per usual: Let 
$i_j: (V_j,||\cdot ||_j) \longrightarrow (X_j,d_j)$ be 
isometries of the normed vector spaces $(V_j,||\cdot ||)$ into the metric spaces $(X_j,d_j)$. Then, with
$||\cdot ||: (V_1\times V_2) \longrightarrow \mathbb{R}$ defined via 
\begin{displaymath}
||(v,w)|| \;:= \; \sqrt{||v||_1^2 \; + \; ||w||_2^2}, \hspace{1.5cm} \forall v\in V_1, w\in V_2,
\end{displaymath} 
the map 
$i:=i_1\times i_2 : (V_1\times V_2, ||\cdot ||)\longrightarrow (X,d):=(X_1\times X_2, \sqrt{d_1^2+d_2^2})$ is an isometry.
Thus $rank_{M}(X_1,X_2) \ge rank_{M} X_1 + rank_{M} X_2$.
\item[ii)] Let $rank_M(X,d)=n$ and let $\varphi : A\longrightarrow X$ be an isometric map, 
where $A$ is the affine space for some $n$-dimensional normed vector space $(V,||\cdot ||)$.
By Proposition \ref{prop} there are two pseudonorms $||\cdot ||_i$, $i=1,2$, on $V$ such 
that $||\cdot ||_1^2+||\cdot ||_2^2=||\cdot ||^2$ and such that 
${\varphi}_i:(A,||\cdot ||_i)\longrightarrow (X_i,d_i)$ are isometric. \\
Let $V_i$ be vector-subspaces transversal to $kern||\cdot ||_i$. Then $dim V_1 + dim V_2\ge n$
and ${\varphi}_i:(V_i,||\cdot ||_i) \longrightarrow X$ are isometric maps. Thus 
$rank_M(X_i,d_i)\ge dimV_i$.
\end{description}
\begin{flushright}
{\bf q.e.d.}
\end{flushright}


\section{Euclidean Rank of Products}

\label{chapter-euclid}

In this section we prove Theorem \ref{theoeukl}. \\
The superadditivity of the Euclidean rank is obvious. Thus it remains to construct an example
such that the equality does not hold. Therefore we construct two norms 
$||\cdot ||_{i}$, $i=1,2$, on $\mathbb{R}^3$, such that 
\begin{description}
\item[i)] there does not exist an isometric embedding of $\mathbb{E}^2$ in 
$(\mathbb{R}^3,||\cdot ||_i)$, $i=1,2$, i.e.,
\begin{displaymath}
rank_E (\mathbb{R}^3,||\cdot ||_i) \; = \; 1, \;\; i=1,2, \hspace{2cm} \mbox{and}
\end{displaymath}
\item[ii)] the diagonal of $(\mathbb{R}^3,||\cdot ||_1)\times (\mathbb{R}^3,||\cdot ||_2)$ is
isometric to the Euclidean space $\mathbb{E}^3=(\mathbb{R}^3,||\cdot ||_e)$, i.e., 
\begin{displaymath}
rank_E\Big( (\mathbb{R}^3,||\cdot ||_1)\times (\mathbb{R}^3,||\cdot ||_2)\Big) \; \ge \; 3.
\end{displaymath}
\end{description}

The norms will be obtained by perturbations of the Euklidean norm $||\cdot ||_e$ in the following way:
\begin{displaymath}
||v||_i \; = \; {\varphi}_i\left( \frac{v}{||v||_e}\right) ||v||_e, \hspace{2cm} 
\forall v\in \mathbb{R}^3,
\end{displaymath}
where the ${\varphi}_i$ are appropriate functions on $S^2$ that satisfy 
${\varphi}_i(\frac{v}{||v||_e})={\varphi}_i(-\frac{v}{||v||_e})$, $i=1,2$, and
${\varphi}_2=\sqrt{2-{\varphi}_1^2}$. Thus their product norm $||\cdot ||_{1,2}$ satisfies
\begin{displaymath}
||(v,v)||_{1,2}^2 \; = \; ||v||_1^2 \; + \; ||v||_2^2 \; = \; {\varphi}^2||v||_e^2 \; + 
\; (2-{\varphi}^2) ||v||_e^2 \; = \; 2 \; ||v||_e^2 
\end{displaymath}
and the diagonal in $(\mathbb{R}^6,||\cdot ||_{1,2})$ is isometric to $\mathbb{E}^3$ and thus
$ii)$ is satisfied. It remains to show that for ${\varphi}_i$ suitable $i)$ holds. \\

Note that for ${\varphi}_i(\frac{v}{||v||_e})=1+{\epsilon}_i (\frac{v}{||v||_e})$, $i=1,2$, with 
${\epsilon}_i$, $D{\epsilon}_i$ and $DD{\epsilon}_i$ sufficiently bounded, the strict convexity of the 
Euclidean unit ball implies strict convexity of the $||\cdot||_i$-unit balls. Since $||\cdot ||_i$ 
is homogeneous by definition it follows that $||\cdot ||_i$, $i=1,2$, are norms. \\

In order to show that $rank_E({\mathbb{R}}^3,||\cdot ||_i)=1$ for suitable functions 
${\varphi}_i=1+{\epsilon}_i$ we
use the following result:
\begin{lemma} \label{ellipse}
Let $(V,||\cdot ||)$ be a normed vector space with strictly convex norm ball and let 
$i:{\mathbb{E}}^2\longrightarrow (V,||\cdot ||)$ be an isometric embedding. Then $i$ is an affine map and
the image of the unit circle in ${\mathbb{E}}^2$ is an ellipse in the affine space $i({\mathbb{E}}^2)$.
\end{lemma}

{\bf \underline{Remark}:} We recall that the notion of an ellipse in a $2$-dimensional vector
space is a notion of affine geometry. It does not depend on a particular norm. Let $A$ be a two 
dimensional affine space on which $V$ acts simply transitively. A subset $W\subset V$ is called an {\it ellipse},
if there are linearly independent vectors $v_1,v_2\in V$ and a point $a\in A$ such that 
\begin{displaymath}
W \; = \; \Big\{ a \; + \; (\cos \alpha \; v_1 \; + \; \sin \alpha \; v_2) \; \Big| \; \alpha \in [0,2\pi ] \Big\} . 
\end{displaymath}

{\bf Proof of Lemma \ref{ellipse}:} \\
In a normed vector space $(V,||\cdot ||)$ the straight lines are geodesics. If the norm ball is strictly 
convex, then these are the unique geodesics. \\
The isometry $i$ maps geodesics onto geodesics and hence straight lines in $\mathbb{E}^2$ onto straight lines
in $V$. Note that the composition of $i$ with an appropriate translation of $V$ yields an isometry that 
maps the origin of $\mathbb{E}^2$ to the origin of $V$. Let us therefore assume that $i$ maps $0$ to 
$0$. It follows that $i$ is homogeneous. Furthermore it is easy to see that parallels are mapped to parallels 
and this finally yields the additivity of $i$ and thus the claim.
\begin{flushright}
{\bf q.e.d.}
\end{flushright}

Now we define functions ${\varphi}_i=1+{\epsilon}_i$ on $S^2$ in a way such that the intersection of the unit 
ball in $({\mathbb{R}}^3,||\cdot ||_i)$  with a $2$-dimensional linear subspace is never an ellipse.
Therefore we will define the ${\epsilon}_i$'s such that their null sets are 8 circles, 4 of 
which are parallel to the equator $\gamma$, the other 4 parallel to a great circle $\delta$
that intersects the equator orthogonally; these null sets being sufficiently close to $\gamma$
and $\delta$ such that each great circle of $S^2$ intersects those circles in at least 8 points.\\
Furthermore no great circle of $S^2$ is completely contained in the null set. \\

Using spherical coordinates $\Theta \in [0,\pi]$, $\Phi \in [0,2\pi]$, $r\in \mathbb{R}^+$,
we define 
\begin{displaymath}
\tilde{\epsilon}_1(\Theta ,\Phi ,r) \; := \; \frac{1}{n} \; \prod\limits_{k=2}^{3} \; 
\sin\Big( \Theta + \frac{k\pi}{8}\Big) \; \sin\Big( \Theta + \frac{(8-k)\pi}{8}\Big),
\end{displaymath}
with $n\in \mathbb{N}$ sufficiently large, such that the norm $||\cdot ||_1$ we will obtain admits 
a strictly convex unit ball. \\
One can easily check that 
$\tilde{\epsilon}_1^k(\Theta ,\Phi ,r)= \sin ( \Theta + \frac{k\pi}{8}) 
\sin (\Theta + \frac{(8-k)\pi}{8})$, $k\in \mathbb{N}$, 
satisfies $\tilde{\epsilon}_1^k(\Theta ,\Phi ,r) = \tilde{\epsilon}_1^k(-\Theta ,\Phi ,r)$ and so
does $\tilde{\epsilon}_1$. Since $\tilde{\epsilon}_1$ is independent of $\Phi$ it satisfies
$\tilde{\epsilon}_1(\frac{v}{||v||})=\tilde{\epsilon}_1(-\frac{v}{||v||})$. Its null set is the union
of the circles parallel to the equator $\gamma$ at 
$\Theta = \{\frac{1}{4}\pi , \frac{3}{8}\pi , \frac{5}{8}\pi , \frac{3}{4}\pi \}$. \\
Define $\hat{\epsilon}_1$ analogous to $\tilde{\epsilon}_1$ but with the null set consisting
of circles parallel to $\delta$ instead of the equator $\gamma$.\\
With that we set ${\varphi}_1 =1+\tilde{\epsilon}_1\hat{\epsilon}_1$ and $||\cdot ||_1$ defined via
\begin{displaymath}
||v||_1 \; := \; {\varphi}_1 \Big( \frac{v}{||v||_e}\Big) \; ||v ||_e
\end{displaymath} 
is a norm on $\mathbb{R}^3$ whose unit ball coincides with the $||\cdot ||_e$-unit ball exactly on the 
null set of $\tilde{\epsilon}_1\hat{\epsilon}_1$. Obviously
\begin{displaymath}
||\cdot ||_2 \; := \; \sqrt{2 \; - \; {\varphi}_1^2 } \; ||\cdot ||_e
\end{displaymath}
is another norm on $\mathbb{R}^3$ whose unit ball also intersects the $||\cdot ||_e$-unit ball 
on the null set of $\tilde{\epsilon}_1\hat{\epsilon}_1$.

\begin{figure}[htbp]
\psfrag{T1}{$\Theta = \frac{\pi}{4}$}
\psfrag{T2}{$\Theta = \frac{3\pi}{8}$}
\psfrag{T3}{$\Theta = \frac{5\pi}{8}$}
\psfrag{T4}{$\Theta = \frac{3\pi}{4}$}
\includegraphics[width=0.9\columnwidth]{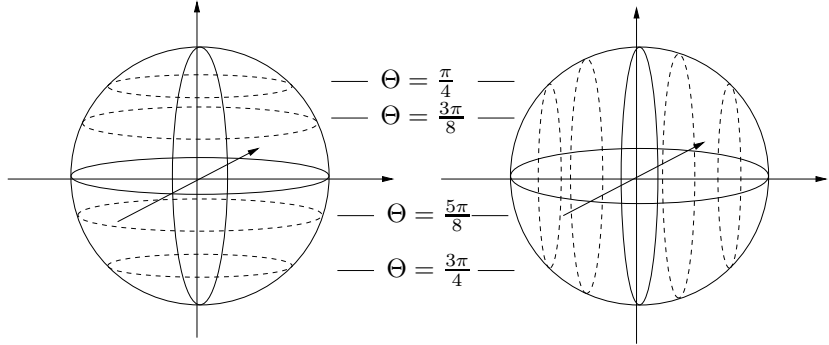}
\caption{The dashed circles in this figure are the sections of the $||\cdot ||_1,||\cdot ||_2$-
and $||\cdot ||_e$ unit balls.}
\end{figure}

We finally conclude that $rank_E({\mathbb{R}}^3,||\cdot ||_j)$, $j=1,2$. \\
Assume to the contrary, that there exists an isometric embedding 
$i:{\mathbb{E}}^2\longrightarrow ({\mathbb{R}}^3,||\cdot ||_j)$. By Lemma \ref{ellipse} we
can assume (after a translation)  that $i$ is a linear isometry and that the image of the unit circle 
$S\subset {\mathbb{E}}^2$ is an ellipse in the linear subspace $i({\mathbb{E}}^2)$ which is in addition
contained in the unit ball $B_j$ of $||\cdot ||_j$. Note that  $i(\mathbb{E}^2)\cap B_j$ and 
$i(\mathbb{E}^2)\cap S^2$ are ellipses which coincide by construction in at least $8$ points.
Since two ellipses with more than $4$ common points coincide, we have 
$i(\mathbb{E}^2)\cap B_j = i(\mathbb{E}^2)\cap S^2$. This contradicts to the fact that by construction
$i(\mathbb{E}^2)\cap B_j \cap S^2$ is a discrete set.
\begin{flushright}
{\bf q.e.d.}
\end{flushright}


\section{Metrics on Product Spaces}
\label{product-metrics}

In section \ref{minkadd} and \ref{chapter-euclid} we discussed the behaviour of the rank
with respect to the usual product. Given a finite number $(X_i,d_i)$, $i=1,...,n$ of metric
spaces there are different possibilities to define a metric $d$ on the product
${\Pi}_{i=1}^{n}X_i$. The standard choice is of course the Euclidean product metric
\begin{displaymath}
d\Big( (x_1,...,x_n),(y_1,...,y_n)\Big) \; = \; 
\Big( \sum\limits_{i=1}^{n} d_i^2(x_i,y_i){\Big)}^{\frac{1}{2}}.
\end{displaymath}
In this section we discuss other natural choices. First of all it is natural to require
that $d\big( (x_1,...,x_n),(y_1,...,y_n)\big)$ only depends on the distances 
$d_i(x_i,y_i)$. \\
We denote by $Q^n:=[0,\infty )^n$ the positive quadrant in 
$\mathbb{R}^n$. On $Q^n$ we define a partial ordering $\le$ in the following way: if 
$q^1=(q_1^1,...,q_n^1)$ and $q^2=(q_1^2,...,q_n^2)$ then
\begin{displaymath}
q^1 \; \le \; q^2 \;\;\; :\Longleftrightarrow \;\;\; q_i^1 \; \le q_i^2 \;\; 
\forall i\in \{ 1,2,...,n\} .
\end{displaymath}
Let $\Phi : Q^n \longrightarrow [0,\infty )$ be a function and consider the function
$d_{\Phi}:X\times X \longrightarrow [0,\infty )$,
\begin{displaymath}
d_{\Phi}\Big( (x_1,...,x_n),(y_1,...,y_n)\Big) \; = \; 
\Phi \Big( d_1(x_1,y_1),...,d_n(x_n,y_n)\Big) .
\end{displaymath} 
In order that  $d_{\Phi}$ will be a metric, we clearly have to assume
\begin{description}
\item[(A)] $\Phi (q) \; \ge \; 0 \;\; \forall q\in Q\;\;\;$ and $\;\;\;\Phi (q) \; = \; 0 \;\; 
\Leftrightarrow \;\; q=0$.
\end{description}
The symmetry of $d_{\Phi}$ is obvious. We now translate the triangle inequality for
$d_{\Phi}$ into a condition on $\Phi$. \\
Let $x=(x_1,...,x_n)$, $y=(y_1,...,y_n)$, $z=(z_1,...,z_n)\in X$ and consider the 
``distance vectors''
\begin{eqnarray*}
q^1 & := & \Big( d_1(x_1,z_1),...,d_n(x_n,z_n)\Big), \\
q^2 & := & \Big( d_1(x_1,y_1),...,d_n(x_n,y_n)\Big) \hspace{1cm} \mbox{and} \\
q^3 & := & \Big( d_1(y_1,z_1),...,d_n(y_n,z_n)\Big)
\end{eqnarray*}
in $Q^n$. Since for every $i\in \{ 1,...,n\}$, $x_i,y_i,z_i$ are points in $X_i$
we see that $q^j\le q^k+q^l$ for every permutation $\{ j,k,l \}$ of $\{ 1,2,3 \}$. \\
Now $d_{\Phi}$ satisfies the triangle inequality if $\Phi$ satisfies
\begin{description}
\item[(B)] for all points $q^1,q^2,q^3\in Q^n$ with $q^j\le q^k+q^l$ we have
\begin{displaymath}
\Phi (q^j) \; \le \;  \Phi (q^k) \; + \; \Phi (q^l).
\end{displaymath}
\end{description}

{\bf \underline{Remark}:}
\begin{description}
\item[i)] Note that for $q^1,q^2,q^3$ one can always take a triple of the form 
$p,q,p+q$, hence (B) implies in particular $\Phi (p+q) \le \Phi (p)+ \Phi (q)$.
\item[ii)] The condition (B) can be applied for the triple $p,q,q$ in the case
that $p\le 2q$. Then $\Phi (p)\le 2 \Phi (q)$. \\
This has the following consequence: If $\Phi$ satisfies (A) and (B) then for every 
$\epsilon >0$ the function $\Phi |_{Q^n\setminus B_{\epsilon}(0)}$ has a positive lower
bound, where $B_{\epsilon}(0)\subset Q^n$ is the $\epsilon$-ball in the 
Euclidean metric. Indeed let $p_i:=\frac{2\epsilon}{\sqrt{n}}e_i$, where $e_i$ is
the unit vector, then for every $q\in Q^n\setminus B_{\epsilon}(0)$ there is $p_i$ with
$2q\ge p_i$ and hence $2\Phi (q) \ge \min \Phi (p_i) > 0$.
\end{description}

It is now easy to prove the following result
\begin{lemma} \label{metric}
Let $\Phi : Q^n\longrightarrow [0,\infty )$ be a function. Then $d_{\Phi}$ is a metric
on $X$ for all possible choices of metric spaces $(X_i,d_i)$, $i=1,...,n$, if 
and only if $\Phi$ satisfies (A) and (B).
\end{lemma}
This Lemma still allows strange metrics on a product (even the trivial
product $n=1$). Let for example $\Phi : Q^n\longrightarrow [0,\infty )$ be an 
arbitrary function with $\Phi (0)=0$ and $\Phi (q)\in\{ 1,2\}$,
$\forall q\in Q^n\setminus \{ 0\}$. Then $d_{\Phi}$ is a metric. \\
If we require however that the product metric space $X$ is always an inner metric space
in the case the $X_i$ are, the conditions on $\Phi$ are very rigid. \\
For convenience of the reader we recall the notion of an inner metric space.
Let $(X,d)$ be a metric space. For a continuous path $c:[0,1]\longrightarrow X$ 
one defines as usual the length
\begin{displaymath}
L(c) \; := \; \sup \Big\{ \sum\limits_{j=1}^k d\Big( c(t_{j-1}),c(t_j) \Big) \Big\},
\end{displaymath}  
where the $sup$ is taken over all subdivisions
\begin{displaymath}
0 \; = \; t_0 \; \le \; t_1 \; \le ... \le \; t_k \; = \; 1 \hspace{1cm} 
\mbox{of} \;\; [0,1].
\end{displaymath}
$(X,d)$ is called an inner metric space if for all $x,y\in X$, $d(x,y)=\inf L(c)$,
where the $\inf$ is taken over all paths from $x$ to $y$.

We need the following 
\begin{lemma} \label{psi-norm}
For $\Phi :Q^n\longrightarrow [0,\infty )$ the function
$\Psi : \mathbb{R}^n \longrightarrow [0,\infty )$ defined via
\begin{displaymath}
\Psi \Big( \sum\limits_{i=1}^n \; x_i \, e_i \Big) \; := \; 
\Phi \Big( \sum\limits_{i=1}^n \; |x_i| \, e_i \Big) 
\end{displaymath}
is a norm on $\mathbb{R}^n$ if and only if $\Phi$ satisfies the following conditions:
\begin{description}
\item[(1)] $\Phi (q)\ge 0$ $\forall q\in Q^n$ and $\Phi (q)=0 \; \Leftrightarrow q=0$,
\item[(2)] $\Phi$ is monoton, i.e. $q\le p \; \Longrightarrow \Phi (q) \le \Phi (q)$
$\forall p,q\in Q^n$,
\item[(3)] $\Phi (p+q) \; \le \; \Phi (p) \; + \; \Phi (q)$,
\item[(4)] $\Phi (\lambda q) \; = \; \lambda \Phi (q)$ $\forall p\in Q^n, \lambda \ge 0$.
\end{description}
\end{lemma}
{\bf Proof of Lemma \ref{psi-norm}:} \\
``$\Longrightarrow$'' \\
Let $\Phi$ satisfy $(1)-(4)$. Then
$\Psi \ge 0$, $\Psi (x)=0 \; \Longleftrightarrow \;x=0$ and 
$\Psi (\lambda x) \; = \; |\lambda | \Psi (x)$ directly follow from the 
definition of $\Psi$. In order to verify the subadditivity, note that
for $x=(x_1,...,x_n)$ and $y=(y_1,...,y_n)$
\begin{eqnarray*}
\Psi (x \; + \; y) & = & \Phi \Big( \sum\limits_{i=1}^n \; |x_i+y_i| \, e_i \Big) \\
& \stackrel{(2)}{\le} & \Phi \Big( \sum\limits_{i=1}^n \; (|x_i|+|y_i|) \, e_i \Big) \\
& \stackrel{(3)}{\le} & \Phi \Big( \sum\limits_{i=1}^n \; |x_i| \, e_i \Big) \; + \; 
\Phi \Big( \sum\limits_{i=1}^n \; |y_1| \, e_i \Big) \\
& =  & \Psi (x) \; + \; \Psi (y).
\end{eqnarray*}
``$\Longleftarrow$'' \\
Assume now that $\Psi$ is a norm. Then $\Phi$ clearly satisfies $(1),(3),(4)$. To prove
(2) it is enough to show that $\Phi (p+\lambda e_i)\ge \Phi (p)$ for any unit vector 
$e_i$ and $\lambda \ge 0$. Assume that $\Phi (p+\lambda e_i)< \Phi (p)$. Write
$p=(p_1,...,p_n)\in Q^n$, let 
$q=(p_1,...,p_{i-1},-p_i-\lambda ,p_{i+1},...,p_n)\in \mathbb{R}^n$. Then
$\Psi (q)=\Psi (p+\lambda e_i)< \Psi (p)$ but $p$ is on the segment between $q$ and 
$p+\lambda e_i$. This contradicts to the subadditivity of $\Psi$. 
\begin{flushright}
{\bf q.e.d.}
\end{flushright}

Now we are able to give the \\
{\bf Proof of Theorem \ref{innermetric}:} \\
``$\Longrightarrow$'' \\
We show that in the case that $(X,d_{\Phi})$ is an inner metric space
for all choices of inner metric spaces $(X_i,d_i)$
$\Phi$ satisfies conditions $(1)-(4)$ as in Lemma \ref{psi-norm}, which then implies
the result. \\

Similar as in Lemma \ref{metric} we see that $\Phi$ satisfies $(A)=(1)$ and $(B)$ which 
implies $(3)$.\\
Consider now the following example: let $(X_i,d_i)=\mathbb{R}$ with the standard metric,
then $X={\mathbb{R}}^n$ and by assumption $d_{\Phi}$ is a length metric on $\mathbb{R}^n$.
The translations of $\mathbb{R}^n$ are isometries of $d_{\Phi}$. $\Phi$ satisfies $(B)$
and hence $\Phi$ assumes a positive lower bound on $Q^n\setminus B_{\epsilon}(0)$
by the remark $ii)$ above. This
implies that a continuous curve $c:[0,1]\longrightarrow (\mathbb{R}^n,d_{\Phi})$ is also
continuous when regarded as a map into $(\mathbb{R}^n,d_e)$, where $d_e$ is the Euclidean
metric. \\

We need the
\begin{sublemma} \label{sublemma1}
Let $q$ be a point in the interior of $Q^n$. For every $\epsilon >0$ there exists a point 
$p\in Q^n$ with 
\begin{displaymath}
d_e\Big( p,\frac{1}{2}q\Big) \; < \; \epsilon \hspace{1cm} \mbox{and} \hspace{1cm}
\Phi (p)\le \frac{1}{2} \Phi (q) \; + \; \epsilon .
\end{displaymath}
\end{sublemma}

{\bf Proof of Sublemma \ref{sublemma1}:}(Claim) \\
Since $(\mathbb{R}^n,d_{\Phi})$ is an inner metric space there is a path 
$c:[0,1]\longrightarrow \mathbb{R}^n$ from $0$ to $q$ continuous in the topology of 
$(\mathbb{R}^n,d_{\Phi})$ and hence also continuous in $(\mathbb{R}^n,d_e)$ with
$L(c)<d(0,q)+\epsilon = \Phi (q) + \epsilon$, where $L$ is the length with respect to 
$d_{\Phi}$. \\
Given a subdivision $0=t_0\le t_1\le ...\le t_k =1$ and a permutation $\pi =(i_1,...,i_k)$
of the numbers $1,...,k$ we define a new continuous path 
$c^{\pi}:[0,1]\longrightarrow \mathbb{R}^n$ from $0$ to $q$ with $L(c^{\pi})=L(c)$ in 
the following way: \\
Let $s_j:=t_j-t_{j-1}$ be the length of the interval $I_j:=[t_{j-1},t_j]$ and let
$m_j=s_{i_1}+...+s_{i_j}$. Define $c^{\pi}$ first on $[0,m_1]=[0,s_{i_1}]$ by
$c^{\pi}(t)=c(t_{i_1-1}+t)-c(t_{i_1-1})$, i.e., we take $c^{\pi}|_{[0,s_{i_1}]}$ as
$c|_{I_{i_1}}$ translated such that $c^{\pi}$ starts at the origin. \\
If $c^{\pi}$ is already defined on $[0,m_{j-1}]$, then define
\begin{displaymath}
c^{\pi}|_{[m_{j-1},m_j]} \; = \; c(t_{i_j-1}+t) \; - \; 
\Big( c^{\pi}(m_{j-1}) \; + \; c(t_{i_j-1})\Big) ,
\end{displaymath}
i.e., $c^{\pi}|_{[m_{j-1}]}$ is $c|_{I_{i_j}}$ translated such that $c^{\pi}$ stays 
continuous. Since translations are isometries one easily sees that $L(c^{\pi})=L(c)$
and clearly $c^{\pi}$ is also a path from $0$ to $q$. \\
It is elementary and not too difficult to prove that there is a subdivision 
$0=t_0\le ...\le t_k=1$ and a permutation $\pi$ such that $c^{\pi}$ stays within
the $\epsilon$-tube (with respect to the Euclidean metric) of the line 
$\mathbb{R}\cdot q\subset {\mathbb{R}}^n$. Thus there is a point $p=c^{\pi}(t_0)$
with $d_e(p,\frac{1}{2}q)<\epsilon$. Let $p'=q-p\in \mathbb{R}^n$,
then $d_e(p',\frac{1}{2}q)=d_e(p,\frac{1}{2}q)<\epsilon$. By choosing 
(in advance) $\epsilon >0$ small enough we can assume that $p'\in Q^n$.  \\
Note that 
\begin{eqnarray*}
\Phi (q) \; + \; \epsilon & \ge & L(c^{\pi}) \\
& = & L(c^{\pi}|_{[0,t_0]}) \; + \; L(c^{\pi}|_{[t_0,1]}) \\
& \ge & \Phi (p) \; + \; \Phi (p') \\
& \ge & \Phi (p+p') \\
& = & \Phi (q).
\end{eqnarray*}
Thus $p$ or $p'$ satisfies the required estimate.
\begin{flushright}
$\Box$
\end{flushright}

\begin{sublemma} \label{sublemma2}
There exists a constant $C>0$ such that $\Phi (v)\le C ||v||_e$, where $||\cdot ||_e$
is the Euclidean norm in $Q^n$. In particular $\Phi $ is continuous at $0$.
\end{sublemma}
{\bf Proof of Sublemma \ref{sublemma2}:} \\
We will show that there exists some constant $C$, such that $\Phi (te_i)\le Ct$ for all
$i=1,...,n$. Then the subadditivity $(4)$ implies the existence of the claimed constant. \\
Thus let $i\in \{ 1,...,n\}$ and $t>0$ be given. Let $C=\Phi \Big( (1,...,1)\Big)>0$.
Choose $k\in \mathbb{N}$ and $m\in \mathbb{N}$ such that 
\begin{equation} \label{*}
\frac{m-1}{2^k} \; \le \; t \; < \; \frac{m}{2^k}.
\end{equation}
By applying Sublemma \ref{sublemma1} several times there exists a point 
$p^k=(p_1^k,...,p_n^k)$ $\epsilon$-close to $(\frac{1}{2^{k+1}},...,\frac{1}{2^{k+1}})$ such
that $\Phi (p^k)\le\frac{C}{2^{k+1}}+\epsilon$, where $\epsilon$ can be chosen as small
as we want. \\
Let $\bar{p}^k=(-p_1^k,...,-p_{i-1}^k,p_i^k,-p_{i+1}^k,...,-p_n^k)\in \mathbb{R}^n$. Then
$p^k+\bar{p}^k=2p_ie_i=:{\rho}_k e_i$, where ${\rho}_k$ is $\epsilon$-close to 
$\frac{1}{2^k}$ and
\begin{eqnarray*}
\Phi ({\rho}_ke_i) & = & d_{\Phi}(0,{\rho}_ke_i) \\
& \le & d_{\Phi}(0,p^k) \; + \; d_{\Phi}(p^k,{\rho}_ke_i) \\
& = & \Phi (p^k) \; + \; d_{\Phi}(0,\bar{p}^k) \\
& = & \Phi (p^k) \; + \; \Phi (p^k).
\end{eqnarray*}
Thus
\begin{equation} \label{**}
\Phi({\rho}_ke_i) \; \le \; \frac{C}{2^k} \; + \; 2\epsilon .
\end{equation}
For $\epsilon$ sufficiently small (depending on $m$) we have 
\begin{displaymath}
(m-2) {\rho}_k \; < \; t \; < \; m {\rho}_k.
\end{displaymath}
Thus we can apply property $(B)$ to the triple 
$\big( (m-2){\rho}_ke_i, te_i, 2{\rho}_ke_i\big)$ and obtain
\begin{eqnarray*}
\Phi (te_i) & \le & \Phi \Big( (m-2){\rho}_ke_i\Big) \; + \; \Phi (2{\rho}_ke_i) \\
& \le & m \Phi ({\rho}_ke_i) \\
& \le & m \Big( \frac{C}{2^k} \; + \; 2\epsilon \Big) \\
& \le & Ct \; + \; \frac{1}{2^k} C \; + \; 2\epsilon m,
\end{eqnarray*}
where the first inequality comes from $(B)$, the second from the subadditivity $(3)$
of $\Phi$, the third from equation (\ref{**}) and the last from equation (\ref{*}). \\
Note that we can make the term $\frac{1}{2^k}C+2\epsilon m$ as small as we like, by 
choosing first $k$ large enough and then choosing $\epsilon$ small enough (depending
on $k$ and $m$). Thus $\Phi(te_i)\le Ct$.
\begin{flushright}
$\Box$
\end{flushright}
The continuity at $0$ and $(B)$ together easily imply that $\Phi$ is continuous
everywhere. \\
Sublemma \ref{sublemma1} and the continuity imply 
$\Phi (\frac{1}{2}q)=\frac{1}{2}\Phi (q)$ $\forall q\in Q^n$. Together with the 
subadditivity $(3)$ this implies $\Phi (\lambda q)=\lambda \Phi (q)$ for all 
dyadic numbers and by continuity for all $\lambda \ge 0$. \\
Finally it remains to show the monotonicity $(2)$. It is enough to show that for
given $p\in Q^n$
\begin{displaymath}
\Phi (p \; + \; \lambda e_i) \; \ge \; \Phi (p) \hspace{1cm} \forall \lambda \ge 0, 
\;\;\; i=1,...,n.
\end{displaymath}
Assume that $\Phi (p+\lambda e_i)<\Phi (p) $ and let 
$\epsilon :=\Phi (p)-\Phi (p+\lambda e_i)>0$. \\
Let us first consider the case that $p_i=0$. Then $2p=(p+\lambda e_i)+(p-\lambda e_i)$,
where $(p-\lambda e_i)\notin Q^n$. Hence 
\begin{eqnarray*}
\Phi (2p) & = & d_{\Phi}(0,2p) \\
& = & d_{\Phi}(0,p \; + \; \lambda e_i) \; + \; d_{\Phi}(0,p \; - \; \lambda e_i) \\
& = & 2\Phi (p \; + \lambda e_i) \\
& > & 2\Phi (p),
\end{eqnarray*}
which is a contradiction to the sublinearity of $\Phi$. \\

If $p_i>0$ then choose $m,k\in \mathbb{N}$ such that 
$0\le mp_i-k\lambda < \frac{\epsilon}{2C}$, where $C$ is the lipschitz constant of 
$\Phi$ with respect to the Euclidean norms on $Q^n$ and $[0,\infty )$. \\
Let $\bar{p}:=\big( p_1,...,p_i,-(p_i+\lambda ),p_{i+1},...,p_n\big)$. Consider now
the points $q:=(2k+2m)p$ and $\bar{q} := (2k+m)(p+\lambda e_i) + m\bar{p} $.
Note that these points only differ at the $i$-th component. The difference at this 
component is
\begin{displaymath}
\Big| (2k\; + \; 2m) p_i \; - \; 2k(p_i \; + \; \lambda )\Big| \; = \; 
2 \Big| mp_i \; - \; k\lambda \Big| \; < \; \frac{\epsilon}{C}.
\end{displaymath}
Thus $|\Phi (q)-\Phi (\bar{q})|<\epsilon$. But $\Phi (q) = (2k+2m)\Phi (p)$ by homogeneity
and
\begin{eqnarray*}
\Phi (\bar{q}) & \le & (2k+m) \Phi (p \; + \; \lambda e_i) \; + \; m d_{\Phi}(0,\bar{p}) \\
& = & (2k+2m) \Phi (p \; + \; \lambda e_i) \\
& < & (2k+2m) \Phi (p) \; - \; \epsilon ;
\end{eqnarray*}
a contradiction. \\

``$\Longleftarrow$'' \\
Let now $\Phi$ satisfy properties $(1)-(4)$. Note that from Lemma \ref{psi-norm} it 
follows immediately that $\Phi$ is continuous with respect to the standard topology
on $Q^n$. In order to show that for any choices of 
inner metric spaces $X_1,...,X_k$ the product $(X,d_{\Phi})$ is an inner metric space we 
prove the following
\begin{lemma} \label{productlength}
Let $(X_i,d_i)$ be metric spaces and $c_i: [0,1]\longrightarrow X_i$ be continuous curves
parametrized by arclength connecting $p_i\in X_i$ with 
$q_i\in X_i$, $i=1,...,k$. Denote by $l_i$ the $(X_i,d_i)$-length of $c_i$ and suppose that
$\Phi$ satisfies conditions $(1)-(4)$. Then the $(X,d_{\Phi})$-length of the product 
curve $c=(c_1,...,c_k):[0,1]\longrightarrow X$ is $L(c)=\Phi(l_1,...,l_k)$. \\
Furthermore $c$ is also parametrized by arclength.
\end{lemma}
{\bf Proof of Lemma \ref{productlength}:} \\
Note that the $(X_i,d_i)$-length $L_{c_i}$ of $c_i$ is given through
\begin{displaymath}
l_i \; = \; L(c_i) \; = \; \lim\limits_{N\longrightarrow \infty}\sum\limits_{j=1}^N \;
d_i \Big( c_i(\frac{j-1}{N}), c_i(\frac{j}{N})\Big) ,
\end{displaymath}
where $d_i\big( c_i(\frac{j-1}{N}), c_i(\frac{j}{N})\big) =\frac{l_i}{N}$. \\
For the $(X,d_{\Phi})$-length $L(c)$ of $c$ one has
\begin{eqnarray*}
L(c) & = & \lim\limits_{N\longrightarrow \infty}\sum\limits_{j=1}^N \;
d_{\Phi} \Big( c(\frac{j-1}{N}), c(\frac{j}{N})\Big) \\
& = &\lim\limits_{N\longrightarrow \infty}\sum\limits_{j=1}^N \; 
\Phi \Big( d_1\Big( c_1(\frac{j-1}{N}),c_1(\frac{j}{N})\Big) ,...,
d_k\Big( c_k(\frac{j-1}{N}),c_k(\frac{j}{N})\Big) \Big) \\
& \stackrel{(2)}{\le} & \lim\limits_{N\longrightarrow \infty}\sum\limits_{j=1}^N \; 
\Phi \Big( \frac{l_1}{N},...,\frac{l_k}{N}\Big) \\
& \stackrel{(4)}{=} & \lim\limits_{N\longrightarrow \infty}\sum\limits_{j=1}^N \; 
\frac{1}{N} \Phi (l_1,...,l_k) \\
& = & \Phi (l_1,...,l_k).
\end{eqnarray*}
On the other hand the continuity and subadditivity of $\Phi$ yield:
\begin{eqnarray*}
L(c) & = & \lim\limits_{N\longrightarrow \infty}\sum\limits_{j=1}^N \;
d_{\Phi} \Big( c(\frac{j-1}{N}), c(\frac{j}{N})\Big) \\
& = &\lim\limits_{N\longrightarrow \infty}\sum\limits_{j=1}^N \; 
\Phi \Big( d_1\Big( c_1(\frac{j-1}{N}),c_1(\frac{j}{N})\Big) ,...,
d_k\Big( c_k(\frac{j-1}{N}),c_k(\frac{j}{N})\Big) \Big) \\
& \stackrel{(3)}{\ge} & \lim\limits_{N\longrightarrow \infty} 
\Phi \Big( \sum\limits_{j=1}^N \; d_1\Big( c_1(\frac{j-1}{N}),c_1(\frac{j}{N})\Big) ,...,
\sum\limits_{j=1}^N \; d_k\Big( c_k(\frac{j-1}{N}),c_k(\frac{j}{N})\Big) \Big) \\
& = & \Phi (l_1,...,l_k),
\end{eqnarray*}
where the last equality is due to the continuity of $\Phi$.
\begin{flushright}
$\Box$
\end{flushright}
Let now $(X_i,d_i)$, $i=1,...,k$, be inner metric spaces. Then the distance of any two points
$p_i,q_i\in (X_i,d_i)$ may be approximated arbitrarily good by the length of continuous 
curves in $(X_i,d_i)$ joining $p_i$ and $q_i$. Thus $(X,d_{\Phi})$ turns out to be
an inner metric space itself, due to the definition of $d_{\Phi}$, the validity of Lemma
\ref{productlength} and the continuity of $\Phi$. \\
This completes the proof of Theorem \ref{innermetric}.
\begin{flushright}
{\bf q.e.d.}
\end{flushright} 

Let us finally focus on products of normed vector spaces and vector spaces with scalar 
products. \\
The analogue of Lemma \ref{psi-norm} in the case of vector spaces with scalar product is
\begin{lemma} \label{psi-sp}
Let $\Phi : Q^n \longrightarrow [0,\infty )$ satisfy $(1)-(4)$ as in Lemma \ref{psi-norm}.
Then $\Psi$ as in Lemma \ref{psi-norm} is induced by a scalar product $g_{\Psi}$
on $\mathbb{R}^n$ if and only if $\Phi$ satisfies the property
\begin{displaymath}
(5) \hspace{1cm}
{\Phi}^2 \Big( \sum\limits_{i=1}^n \; {\lambda}_i \, e_i\Big) \; = \; 
\sum\limits_{i=1}^{n} \; {\Phi}^2 ({\lambda}_i \, e_i) \hspace{1cm} \forall {\lambda}_i >0.
\end{displaymath}
In this case the set $\{ e_1,...,e_n \}$ is an orthogonal system of $g_{\Psi}$.
\end{lemma}
{\bf  Proof of Lemma \ref{psi-sp}:} \\
From Lemma \ref{psi-norm} we know that $\Psi$ is a norm if and only if
conditions $(1)-(4)$ hold. Now we show that $\Psi$ satisfies the parallelogram
equation if and only if condition $(5)$ also holds: \\
``$\Longrightarrow$'' \\
Suppose that $\Phi$ satisfies condition $(5)$. Then for $x=(x_1,...,x_n)$ and 
$y=(y_1,...,y_n)$ the parallelogram equation is equivalent to
\begin{displaymath}
\sum\limits_{i=1}^n \; \Big[ |x_i+y_i|^2 \; + \; |x_i-y_i|^2 \; - \; 
2 \, \big[ |x_i|^2 \, + \, |y_i|^2 \big] \Big] {\Phi}^2 (e_i) \; = \; 0 ,
\end{displaymath} 
which holds trivially. \\
``$\Longleftarrow$'' \\
Now suppose that the parallelogram equation holds. For $x=(x_1,...,x_{n-1},0)$ and
$y=(0,...,0,y_n)$ it takes the form 
\begin{displaymath}
{\Phi}^2 \Big( \sum\limits_{i=1}^{n-1} \; |x_i| \, e_i \; + \; |y_n| \, e_n \Big) \; = \;
{\Phi}^2 \Big( \sum\limits_{i=1}^{n-1} \; |x_i| \, e_i \Big) \; + \; 
{\Phi}^2 \Big( |y_n| \, e_n \Big) .
\end{displaymath}
The same computation for $x=(x_1,...,x_{n-2},0,0)$, $y=(0,...,0,y_{n-1},0)$ and so on finally 
yields condition $(5)$.
\begin{flushright}
{\bf q.e.d.}
\end{flushright}

Since a normed vector space with the metric induced by its norm is an inner metric space,
we can use Theorem \ref{innermetric} in order to prove the following results:
\begin{description}
\item[i)] Let $(V_i,||\cdot ||_i)$, $i=1,...,k$ be normed vector spaces and 
$\Phi : Q^n\longrightarrow [0,\infty )$ be a function. Define the function
$||\cdot ||_{\Phi}:V=V_1\times ... \times V_k \longrightarrow [0,\infty )$ through
\begin{displaymath}
\Big| \Big| (v_1,...,v_k)\Big| {\Big|}_{\Phi} \; := \;  \Phi 
\Big( \sum\limits_{i=1}^k \; ||v_i||_i \, e_i \Big) .
\end{displaymath} 
Then $(V,||\cdot ||_{\Phi})$ is a normed vector space for all possible choices of normed 
vector spaces $(V_i,||\cdot ||_i)$ if and only if $\Psi$ as defined in Lemma 
\ref{psi-norm} is a norm.
\item[ii)] Let $(V_i,||\cdot ||_i)$, $i=1,...,k$, be normed vector spaces the norms of which 
are induced by scalar products $<\cdot ,\cdot >_i$ on $V_i$ and
$\Phi : Q^n\longrightarrow [0,\infty )$ be a function. \\
Then the norm $||\cdot ||_{\Phi}$ on $V=V_1\times ... \times V_k$ as in $i)$ is 
induced by a scalar product $<\cdot ,\cdot >_{\Phi}$ for all choices
of vector spaces $V_i$ with scalar products $<\cdot ,\cdot >_i$, if and only if
the norm $\Psi$ as defined in Lemma \ref{psi-norm} is induced by a scalar product
$g_{\Psi}$ on $\mathbb{R}^n$. \\
Thus for two vectors $v=(v_1,...,v_k),w=(w_1,...w_k)\in v$ one always has
\begin{displaymath}
<v,w>_{\Phi} \; = \; \sum\limits_{i=1}^k \; {\Phi}^2 (e_i) \; <v_i,w_i> \; ,
\end{displaymath}
which is the usual Euclidean product up to a scale of the scalar products on the factors. \\
Note that the degree to that $\{ e_1,...,e_n\}$ fails to be an orthonormal basis
of $g_{\Psi}$ is the degree to that $<\cdot ,\cdot >_{\Phi}$ differs from the
standard scalar product of Euclidean products.
\end{description}


\section{Minkowski Rank of Products II}

In section \ref{minkadd} we proved the additivity of the Minkowski rank with respect to
standard metric products. In this section we want to generalize this theorem to the case
of more general metric products. As the proof
is almost the same as the one of Theorem \ref{theomink} we are only going to comment on 
those parts of the proof that involve new aspects. \\

The analogue of Proposition \ref{prop} is
\begin{proposition} \label{generalized-prop}
Let $A$ denote an affine space on which the normed vector space $(V,|\cdot |)$ acts
simply transitively. Let further $(X_i,d_i)$, $i=1,...,n$, be metric spaces, 
$\Phi :Q^n \longrightarrow [0,\infty )$ be a function satisfying
conditions $(1)-(4)$ such that the norm ball of $\Psi$ is strictly convex and 
let $\varphi : (A,|\cdot |) \longrightarrow ({\Pi}_{i=1}^nX_i,d_{\Phi})$ be an 
isometric map. Then there exist pseudonorms $||\cdot ||_i$, $i=1,...,n$, on V such
that
\begin{description}
\item[i)] $|v| \; = \; \Phi \Big( \sum\limits_{i=1}^n \; ||v||_i \, e_i \Big)
\hspace{1cm} \forall v\in V \hspace{1cm} \mbox{and}$
\item[ii)] ${\varphi}_i:(A,||\cdot ||_i)\longrightarrow (X_i,d_i)$, $i=1,...,n$ are
isometric.
\end{description}
\end{proposition} 

Once again we define ${\alpha}_i:A\times V \longrightarrow [0,\infty)$, $i=1,...,n$, via
\begin{displaymath}
{\alpha}_i(a,v) \; := \; d_i\Big( {\varphi}_i(a),{\varphi}_i(a+v) \Big).
\end{displaymath}
Since $\varphi$ is isometric we have
\begin{equation} \label{lemma1'1}
\Phi \Big( \sum\limits_{i=1}^n \; {\alpha}_i(a,v) \, e_i \Big) \; := \; 
d_{\Phi} \Big( \varphi (a),\varphi (a+v)\Big) \; = \; |v|.
\end{equation}
In order to prove Proposition \ref{generalized-prop} one proves four lemmata,
Lemmata $\ref{lemma1}',\ref{lemma2}',\ref{lemma3}'$ and $\ref{lemma4}'$, that read exactly
like their analogous in the standard product case 
(Lemmata \ref{lemma1}, \ref{lemma2}, \ref{lemma3} and \ref{lemma4} ). \\
While the proofs of Lemmata $\ref{lemma2}',\ref{lemma3}'$ and $\ref{lemma4}'$ are
straight forward to rewrite from those of Lemmata 
\ref{lemma2}, \ref{lemma3} and \ref{lemma4}, we are going to give the proof of 
Lemma $\ref{lemma1}'$, as this involves the extra condition that the norm ball
of $\Psi$ is strictly convex: \\
{\bf Proof of Lemma $\ref{lemma1}'$:} \\
The $d_i$'s triangle inequality yields
\begin{displaymath}
{\alpha}_i(a,v) \; + \; {\alpha}_i(a+v,v) \; \ge \; {\alpha}_i(a,2v). 
\end{displaymath}
Therefore the monotonicity of $\Phi$ gives
\begin{equation} \label{lemma1'2}
\Phi \Big( \sum\limits_{i=1}^n \; [{\alpha}_i(a,v) \; + \; {\alpha}_i(a+v,v) ] \, e_i \Big)
\; \ge \; 
\Phi \Big( \sum\limits_{i=1}^n \; {\alpha}_i(a,2v) \, e_i \Big)
\; = \; 
|2v| \; = \; 2 \, |v|.
\end{equation}
Set $x:=\sum_{i=1}^n{\alpha}_i(a,v)e_i$ and $y:=\sum_{i=1}^n{\alpha}_i(a+v,v)e_i$ and
note that with equations (\ref{lemma1'1}) and (\ref{lemma1'2}) one has
\begin{displaymath}
\Phi (x) \; = \; |v| \; = \; \Phi (y) \hspace{1cm} \mbox{and} \hspace{1cm}
\Phi (x+y) \; = \; 2 \; |v| 
\end{displaymath}
and hence
\begin{displaymath}
\Phi (x+y) \; = \; \Phi (x) \; + \; \Phi (y) \; = \; 2\Phi (x) \; = \; 2\Phi (y).
\end{displaymath}
From this it follows with the strict convexity of $\Phi$ that
\begin{displaymath}
x \; = \; y \hspace{0.5cm} \Longleftrightarrow \hspace{0.5cm} 
{\alpha}_i(a,v) \; = \; {\alpha}_i(a+v,v) \;\;\; \forall i=1,...,n,
\end{displaymath}
which proves Lemma $\ref{lemma1}'$.
\begin{flushright}
{\bf q.e.d.}
\end{flushright}
Using Proposition \ref{generalized-prop} one can finally prove Theorem 
\ref{generalized-minkadd}. 
The proof, however, is just along the lines of the proof of Theorem \ref{theomink} and
will therefore be omitted here.

{\footnotesize UNIVERSIT\"AT Z\"URICH, MATHEMATISCHES INSTITUT, WINTERTHURERSTRASSE 190, 
CH-8057 Z\"URICH, SWITZERLAND \\
E-mail addresses: $\;\;\;\;\;$ foertsch@math.unizh.ch $\;\;\;\;\;$ vschroed@math.unizh.ch}

\end{document}